\documentclass[12pt]{article}
\usepackage[dvips]{graphicx}
\usepackage{amssymb,amsmath}

\newtheorem{thm}{Theorem}[section]
\newtheorem{prop}[thm]{Proposition}
\newtheorem{cor}[thm]{Corollary}
\newtheorem{lem}[thm]{Lemma}
\newtheorem{conj}[thm]{Conjecture}
\newtheorem{exa}[thm]{Example}
\newtheorem{claim}[thm]{Claim}
\newtheorem{notat}[thm]{Notation}
\newtheorem{alg}[thm]{Algorithm}

\newtheorem{defn}[thm]{Definition}
\newtheorem{rem}[thm]{Remark}

\newbox\proofbox
\def\proof{\futurelet\next\lookforbracket}
\def\lookforbracket{\ifx\next[\let\go\usespecialterm
\else\let\go\relax
\ifvmode\vskip-\lastskip\fi\global\setbox\proofbox=\vbox
\bgroup%
\noindent{\\\sc Proof:}%
\enskip\relax\fi\ignorespaces\go
}\def\usespecialterm[#1]{\ifvmode\vskip-\lastskip\fi
\global\setbox\proofbox=\vbox\bgroup%
\noindent\hskip\parindent%
{\\\noindent\sc Proof {#1}:}\textnormal\ \ \relax
\ignorespaces}
\def\endproof{\hbadness10000\parfillskip=0pt\egroup
\unvbox\proofbox
\setbox0=\lastbox
\ifdim\ht0>1pt 
\vskip-2pt
\noindent
\hbox to\textwidth{\vbox{
\parindent=0pt
\unhbox0 $\square$ \hss}\hss}
\relax%

\fi}

\newcommand{\ben}{\begin{enumerate}}
\newcommand{\een}{\end{enumerate}}
\newcommand{\ble}{\begin{lem}}
\newcommand{\ele}{\end{lem}}
\newcommand{\bt}{\begin{thm}}
\newcommand{\et}{\end{thm}}
\newcommand{\bpr}{\begin{prop}}
\newcommand{\epr}{\end{prop}}
\newcommand{\bco}{\begin{cor}}
\newcommand{\eco}{\end{cor}}
\newcommand{\bcon}{\begin{conj}}
\newcommand{\econ}{\end{conj}}
\newcommand{\bde}{\begin{defn}}
\newcommand{\ede}{\end{defn}}
\newcommand{\bex}{\begin{exa}}
\newcommand{\eex}{\end{exa}}
\newcommand{\brem}{\begin{rem}}
\newcommand{\erem}{\end{rem}}
\newcommand{\bnot}{\begin{notat}}
\newcommand{\enot}{\end{notat}}
\newcommand{\balg}{\begin{alg}}
\newcommand{\ealg}{\end{alg}}
\newcommand{\bcl}{\begin{claim}}
\newcommand{\ecl}{\end{claim}}

\renewcommand{\>}{\right > }

\newcommand{\N}{\mathbb{N}}
\newcommand{\C}{\mathbb{C}}
\newcommand{\PP}{\mathbb{P}}
\newcommand{\Z}{\mathbb{Z}}

\begin{document}
\title{Conjugacy in Artin Groups and Application to the 
Classification of Surfaces}
\author{Godelle Eddy, Kaplan Shmuel\footnote{Partially supported by the Emmy Noether Research Institute for
Mathematics, Bar-Ilan University and the Minerva Foundation, Germany, and by the Excellency Center ``Group 
theoretic methods in the study of algebraic varieties'' of the National Science Foundation of Israel.}
\footnote{The paper is part of the author's PhD. Thesis} and Teicher Mina}
\date{}
\maketitle

%%%%%%%%%%%%%%%%%%%%%%%%%%%%%%%%%%%%%%%%%%%%%%%%%%%%%%%%%%%%%%%%%%%%%%%%%%%%%%%%%%%%%%%%%%%%%%%%%%%%%%%%%%%%%
%%%%%%%%%%%%
\begin{abstract} We show that the double reversing algorithm proposed
by Dehornoy in \cite{Deh2} for solving the word problem in
the braid group can also be used to recognize the conjugates of
powers of the generators in an Artin group of spherical type. The proof uses a
characterization of these powers in terms of their
fractional decomposition. This algorithm could have potential applications to braid-based cryptography; it also provides a fast method for testing a necessary condition in the classification of surfaces in algebraic geometry.
\end{abstract}

%%%%%%%%%%%%%%%%%%%%%%%%%%%%%%%%%%%%%%%%%%%%%%%%%%%%%%%%%%%%%%%%%%%%%%%%%%%%%%%%%%%%%%%%%%%%%%%%%%%%%%%%%%%%%
%%%%%%%%%%%%
\section{Introduction and motivation}

Let $S$ be a finite set and let $M=(m_{s,t})_{s,t \in S}$ be a Coxeter matrix, that is, a symmetric matrix 
for which $m_{s,s}=1$ for $s$ in $S$ and $m_{s,t}$  $\N \cup\{ \infty\} \setminus \{0,1\}$ where $s \neq t$. 
The \emph{Artin-Tits system associated to $M$} is the pair $(A,S)$ where $A$ is the group with presentation:

\begin{equation} \label{present} A=\left < S \mid \underbrace{s\ t\ s\ \cdots}_{m_{s,t} \ letters} = 
\underbrace{t\ s\ t\ \cdots}_{m_{s,t} \ letters} \ for \ s,t\textrm{ in } S \textrm{ with } s\neq t \textrm{ 
and } m_{s,t} \neq \infty \right >\end{equation}
$A$ is called the \emph{Artin group associated to $M$}. For
instance, if
$S = \{s_1,\cdots, s_n\}$ with $m_{s_i,s_j}= 3$ for $|i-j| =1$ and
$m_{s_i,s_j} = 2$ otherwise, then the associated Artin group is the braid group on $n+1$ strings, $B_{n+1}$.
If we add the relations $s^2=1$ for every $s$ in $S$ to the presentation of $A$, we get the Coxeter group 
$W_S$ associated with $A$. When $W_S$ is finite, we say that $A$ is of \emph{spherical type}. For instance, 
the braid group is of spherical type. The \emph{Artin monoid~$A^+$ associated to $M$} is the sub-monoid of 
$A$ generated by $S$. This monoid has the same presentation as $A$ but as monoid (see \cite{Paris2}). When 
$A$ is of spherical type,
any element $\beta$ of $A$ can be written in a unique way as an irreducible  fraction $\beta = ab^{-1}$ with 
$a$ and $b$ in $A^+$ (see next section for precise definition). This fraction is called the \emph{orthogonal 
splitting of $\beta$}.  

 The main result of this paper is the following:

\bt \label{THcentral}
Let $(A,S)$ be an Artin-Tits system where $A$ is of spherical type. If $\beta$ is in $A$ and is conjugated to a 
power of some element of $S$, then the orthogonal splitting of $\beta$ is $\beta=bs^jb^{-1}$ with $b$ in 
$A^+$,   $s$ in $S$ and $j$ in $\Z$.
\et
We will derive from this theorem a quadratic (on the length of $\beta$) algorithm which detects these 
conjugated elements and returns $b$, $s$ and $j$ as above. Finding a fast algorithm which recognizes the conjugates of
powers of the generators is  useful both from a theoretic point of view (see \cite{God_Paris} for instance) and from a practical point of view since the crypto-systems which use the braid group are based on the conjugacy problem. Also, following an idea suggested in \cite{BGTI} and in \cite{BGTII}, it was proved in  \cite{KUL} that 
\emph{braid monodromy factorization} can be used as a tool to classify surfaces up to diffeomorphism (see 
last section for more details). As an application of our current result, we deduce a second algorithm  testing for a factorization 
to be a braid monodromy factorization of a cuspidal curve.

The paper is organized as follows: in Section $2$ we prove Theorem $1.1$, in Section $3$, we explain the quadratic algorithm and compute its complexity. Finally, Section $4$ is devoted to the application.

%%%%%%%%%%%%%%%%%%%%%%%%%%%%%%%%%%%%%%%%%%%%%%%%%%%%%%%%%%%%%%%%%%%%%%%%%%%%%%%%%%%%%%%%%%%%%%%%%%%

\section{Proof of the Theorem \ref{THcentral}}

\subsection{Preliminaries}
In this part we recall some basic results on Artin groups which are needed in the proof of Theorem 
\ref{THcentral}.

Let $(A,S)$ be an Artin-Tits system where $A$ is of spherical type, and
let $A^+$ be the associated monoid of $A$. As both words in any relation of the presentation of $A^+$ have 
the same number of letters,  we can define on $A^+$
a unique length function $\ell: A^+\to \mathbb{N}$ which is a morphism of monoids and such that the length of 
an element of $S$ is $1$. We say that $\alpha$ \emph{right-divides} $\beta$ in $A^+$ if $\beta = 
\gamma\alpha$ for some $\gamma$ in $A^+$ ; in that case, we write $\beta\succ\alpha$. We can define in the 
same way the left-divisibility and write $\alpha\prec\beta$ when $\alpha$ \emph{left-divides} $\beta$ in 
$A^+$. Recall also that a monoid $G^+$ is said to be cancellative if any equality  $\beta\alpha_1\gamma = 
\beta\alpha_2\gamma$ in $G^+$ implies $\alpha_1 = \alpha_2$.

\bpr[\cite{BrS}]\label{prop1} 
Let $(A,S)$ be an Artin-Tits system where $A$ is of spherical type and $A^+$ be the associated monoid of 
$A$.\\
i) $A^+$ is cancellative.\\
ii) any two elements of $A^+$ have a right-\emph{l.c.m.} $\alpha\lor_\succ\beta$ and a \emph{right-g.c.d.} 
$\alpha\land_\succ\beta$  in $A^+$ for right-divisibility.\\ 
The same is true if we replace right by left.\\
iii) The right-gcd and the left-gcd of all elements in $S$ are equal. We denoted by $\Delta_S$ this 
element.\\
iv) if $\alpha s = \beta t$ in $A^+$ with $s$ and $t$ in $S$ distinct, then $m_{s,t}\neq\infty$ and there 
exists
$\gamma$ in $A^+$ such that $\alpha  =\gamma\underbrace{\cdots\ t\ s\ t}_{m_{s,t} -1\textrm{ letters}}$ and 
$\beta =
\gamma\underbrace{\cdots\ s\ t\ s}_{m_{s,t} -1\textrm{ letters}}$\\
v) $\Delta_S$ is quasi-central in $A_S$ and $\Delta_S^2$ is in the center in $A_S$; in other words, there 
exists a permutation  $\sigma : S\to S$ of order $1$ or $2$ such that  $s\Delta_S = \Delta_S \sigma(s)$ for any $s$ of 
$S$. 
\epr

The key property  of Artin groups of spherical type is the existence of the element $\Delta_S$ which allows 
to describe the group using the monoid:
\bpr \label{prop2} Let $(A,S)$ be an Artin-Tits system where $A$ is of spherical type and let $A^+$ be the 
associated monoid of $A$.\\
i) (\cite{Del} Paragraph 4) Let $\beta$ be in $A$; then $\beta$ can be written
 $\beta = \beta_1\Delta_S^{-n}$ with $\beta_1$ in $A^+$, and $n$ in $\mathbb{N}$.\\
ii) (\cite{Cha1} Theorem 2.6) Let
 $\beta$ be in $A$; then there exists a unique pair $(a,b)$ of elements of $A^+$ such that $a\land_\succ b 
= 1$ and $\beta = ab^{-1}$. \epr
We call the decomposition  $\beta = ab^{-1}$ of ii) the (right) orthogonal
splitting of $\beta$. In a similar way one can define the left orthogonal splitting $\beta = c^{-1}d$ of 
$\beta$.

Of course, since $\Delta_S$ is quasi-central, we can also in i) write  $\beta = \Delta_S^{-n}\beta_2$ with 
$\beta_2$ in $A^+$.\\

With this background, we are now ready to prove Theorem \ref{THcentral}.
\subsection{Proof of Theorem \ref{THcentral}} We  assume in the following that $A_S$ is an Artin group of 
spherical type, and we choose
$\beta$ in $A_S$ which is conjugated to a power of some
element of $S$, that is $\beta = wt^jw^{-1}$ with $w$ in $A_S$, $t$ in $S$ and $j$ in $\mathbb{Z}$.\\
We will use the following easy result.
\ble \label{lemmeBrS}let $t$ and $u$ be distinct in $S$, let $\alpha$ be in $A^+$ and $k$ in 
$\mathbb{N}-\{0\}$  such that $u$ right-divides $\alpha t^k$. Then $m_{t,u}$ is finite and 
$\underbrace{\cdots\ u\ t\ u}_{m_{t,u}-1\ letters}$ right-divides $\alpha$.  
\ele  

\begin{proof}
This lemma is a direct consequence of Lemma 3.1 and Lemma 3.2 of \cite{BrS} and is true in any, not 
necessarily spherical type Artin monoid. It can also be seen as a consequence of Proposition \ref{prop1} i) 
and iv).
\end{proof}

\begin{proof}(Theorem \ref{THcentral})
If $ab^{-1}$ is the orthogonal splitting of $\beta$ then  the orthogonal splitting of $\beta^{-1} = 
wt^{-j}w^{-1}$ is $ba^{-1}$. Hence it is enough to prove Theorem \ref{THcentral} when $j$ is positive (the 
case $j = 0$ is obvious). So we assume $j$ positive in the following.\\
By  Proposition \ref{prop2} i), there exists
$w_1$ in $A^+$ such that $w = w_1\Delta_S^{-k}$ with $k$ in $\mathbb{N}$. Furthermore, Proposition 
\ref{prop1} vi) states that $\Delta_S$ is quasi-central. Thus  $\beta = wt^jw^{-1}  =
w_1\Delta_S^{-k}t^j\Delta_S^kw_1^{-1} =
w_1\left(\sigma^k(t)\right)^jw_1^{-1}$ with $\sigma^k(t)$ which is in
$S$. Hence we can assume for the remaining of the proof that $w = w_1$
is in $A^+$. Let us prove that the orthogonal splitting of $\beta$ is $bs^jb^{-1}$ for some $b$ in $A^+$ and 
some $s$ in $S$ by induction on $\ell(w)$. If $wt^j\land_\succ w = 1$ then  $wt^jw^{-1}$ is the orthogonal 
splitting of $\beta$ and we are done.  This is, in particular, the case when $w = 1$; so assume $\ell(w)\geq 
1$ and $wt^j\land_\succ w \neq 1$. Assume that the orthogonal splitting of $w's^j{w'}^{-1}$  for any $w'$ in 
$A^+$ with $\ell(w')<\ell(w)$ and any $s$ in $S$ is of the required form. Since $wt^j\land_\succ w$ is 
different from $1$, there exists $u$ in $S$ and $w'$ in $A^+$ such that $w = w'u$ and $w't^j \succ u$. If $u 
= t$ then $\beta = w't^j{w'}^{-1}$ and we can apply the induction hypothesis  since $\ell(w') = \ell(w)-1$. 
Assume that $u$ is different from $t$, then by Lemma \ref{lemmeBrS} we have $w = w''\underbrace{\cdots\  u\ 
t\ u}_{m_{t,u}-1\ letters}$ and $\beta = wt^jw^{-1}  = w''\underbrace{\cdots\ u\ t\ u}_{m_{t,u}-1\ 
letters}t^j\underbrace{u^{-1}t^{-1}u^{-1}\cdots}_{m_{t,u}-1\ letters} {w''}^{-1} = w''s^j{w''}^{-1}$ with $s$ 
in $\{u,t\}$ such that $\underbrace{\cdots\ t\ u\ t}_{m_{t,u}\ letters} = s\underbrace{\cdots\ u\ t\ 
u}_{m_{t,u}-1\ letters}$. Since $\ell(w'') = \ell(w)- m_{t,u} + 1 < \ell(w)$, we can apply the induction 
hypothesis to conclude.
\end{proof}

Note that we can not expect to extend Theorem \ref{THcentral} to the family of Garside groups, which is a 
larger family than the one of Artin groups of spherical type (see \cite{Deh1}), since in $\left< a,b \mid aba 
= b^2\right>$ the normal form of $ab^2a^{-1}$ is $a^2b$.
%%%%%%%%%%%%%%%%%%%%%%%%%%%%%%%%%%%%%%%%%%%%%%%%%%%%%%%%%%%%%%%%%%%%%%%%%%%%%%%%%%%%%%%%%%%%%%%%%%%
\section{The algorithm}

Before we describe our algorithm we note that an algorithm for solving the conjugacy problem in 
Garside groups was already given in \cite{meneses}. Their solution, although not polynomial in general is 
polynomial in the case of the conjugacy class of powers of the group's generators. We present here a 
different approach, which is at least as efficient.

\subsection{Description of the algorithm}

We begin by recalling from \cite{Deh2} and \cite{Deh1} the details of the double reversing method proposed  by Dehornoy to find the orthogonal 
splitting of $\beta$.
    
An Artin group of spherical type is a Garside group with $\Delta_S$ as it's Garside element(\cite{Deh1}). In 
particular, $M_S$ the set of all left divisors of $\Delta_S$ is closed under right-lcm and under left-lcm, 
and includes $S$; the permutation $\sigma$ of Proposition \ref{prop1} v) can be extended to a permutation 
$\sigma$ with order 2 of $M_S$ such that for any $x$ in $M_S$, one has $\sigma(x)\Delta_S = \Delta_S x$. 
Furthermore, there exists a (right) \emph{complement} function $f : M_S \times M_S \to M_S$ such that 
$f(x,x)$ is the empty word for every $x$ in $M_S$ and $xf(y,x)=yf(x,y) = x\land_\prec y$ for any $x,y$ in 
$M_S$.

One has:

\bpr Let $(A,S)$ be an Artin-Tits system where  $A$ is of spherical type; Then, $A$ admits the following 
presentation:
$$\<M_S | xf(y,x)=yf(x,y) \text{ for all } x,y \in M_S \>$$
and $A^+$ admits the same presentation as a monoid.\epr

In the following, we fix $(A_S,S)$ an Artin-Tits system where $A_S$ is of spherical type. If $X$ is a finite 
set, we denote  by $F^*(X)$  the free monoid with base $X$. For $w$ in $F^*(M_S)$ or in $F^*(M_S \cup 
M_S^{-1})$  we write $\overline w$ to denote it's corresponding element in $A$. The elements of these monoids are 
called words; we denote by $\epsilon$ the empty word.  If $\beta$ is in $A$ and $w$ is in $F^*(M_S \cup 
M_S^{-1})$, we say that $w$ represents $\beta$ when $\beta = \overline w$. The complexity of our algorithm is 
calculated regarding the length $\ell$ on $F^*(M_S \cup M_S^{-1})$. In particular, we assume that the 
$f$-reversing function is known. Note that  $S$ is a subset of $M_S$. 
\bde
Let $w=x_1^{\epsilon _1} \cdots x_l^{\epsilon _l}$ be in  $F^*(M_S \cup M_S^{-1})$  such that the $x_i$ are 
in $M_S$ and $\epsilon_i$ are in $\{-1,1\}$. Let $w'=ab^{-1}$  be in  $F^*(M_S \cup M_S^{-1})$  with $a$ and 
$b$ in $F^*(M_S)$; hence $w'$ does not contain any pair $x^{-1}y$  with $x,y$ in $M_S$. We say that $w$ is 
\emph{$f$-reversible} into $w'$ if there exists a finite sequence of $w=w_0,w_1,\cdots, w_n=w'$ such that in 
step $i$ some subword $x^{-1}y$  with $x,y$ in $M_S$ of $w_{i-1}$ is replaced by $f(y,x)f(x,y)^{-1}$ in order 
to get $w_i$.
\ede

\brem \label{THcomplex}
The process of $f$-reversing (on the right) can be duplicated into $g$-reversing on the left using a left 
complement function. In this case, during step $i$ a subword $xy^{-1}$  with $x,y$ in $M_S$ of $w_{i-1}$ is 
replaced by $g(y,x)^{-1}g(x,y)$ in order to get $w_i$, hence in $w'$ there will be no pairs $xy^{-1}$ with 
$x,y$ in $M_S$. Note that for an Artin group of spherical type, it is the same to know the $f$-reversing 
function and the $g$-reversing function using the symmetry of the relations.
\erem

It is well known that every word $w$, that represents $\beta$ in $A$, is $f$-reversible on the right and 
$g$-reversible on the left. Therefore, given $\beta$ in $A$ and $z$ that represents $\beta$, one can use the $g$-reversing process on 
$z$ in order to write $w=w_1^{-1}w_2$ with $w_1,w_2$ in $F^*(M_S)$ which represents $\beta$, and then activate 
the $f$-reversing process on $w$ resulting with $ab^{-1}$. By \cite{Deh3} Proposition 2.12 $\overline{a}{\overline{b}}^{-1}$ is the orthogonal splitting of $\beta$. 

So our algorithm to identify and compute the root of a generator in an Artin group of spherical type is as 
follows:

\balg \label{ALG} (Identify generator conjugation)\\
{\bf Input:} An element $\beta$ of $A$ represented by $w_0=s_1^{\epsilon _1} \cdots s_{\ell}^{\epsilon 
_\ell}\in A$ where the $s_i$ are in $S$ and the $\epsilon_i$ are in $\{-1,1\}$.\\
{\bf Output:} True, if $\beta=bs^jb^{-1}$ (for some element $b$ in $A$ and some $j$), $b$, $s$ and $j$; False 
otherwise.
\parskip=0pt
\ben
\item
Write $w=ab^{-1}$ using the process described above such that $\overline{a}\overline{b}^{-1}$ is the 
orthogonal splitting of $\beta$.
\item
Repeat step $(2)$ on $b^{-1}a$, getting its orthogonal splitting $a_1b_1^{-1}$
\item
If $b_1= \epsilon$ and $a_1=s^j$ for some $s$ in $S$ and $j$ in $\mathbb{N}$. Return \emph{True}, $j,b$ and 
$s$;\\
If $a_1 = \epsilon$  and $b_1=s^j$ for some $s$ in $S$ and $j$ in $\mathbb{N}$. Return \emph{True}, $-j,a$ 
and $s$;\\Else return \emph{False}.
\een
\ealg

\subsection{Complexity}

\bpr \label{prop4}
The process of finding the orthogonal splitting of $\beta$ in $A$ represents by the word $w = 
x_1^{\epsilon_1} \cdots x_\ell^{\epsilon_\ell}$ with $x_i$ in $M_S$ is bounded by $N\ell^2/2$ steps, where 
$N$ is an upper bound to the number of steps needed to compute the completion function $f(x,y)$. Moreover, 
the length (on $M_S$) of the word which represents $\beta$ is not greater than $\ell$.
\epr

\begin{proof}
Since the relations in an Artin group are symmetrical the upper bound is the same for both $f(x,y)$ and 
$g(x,y)$. Now, the First step in the process of computing the orthogonal splitting of $\beta$, is to write 
$\beta$ as $w_1^{-1}w_2$ using the $g$-reversing process on the left. For Second step one activates the 
$f$-reversing process on the result of First step. Because in each step of the $f$-reversing process (and the 
$g$-reversing process) a subword $x^{-1}y$ ({\it resp.} $xy^{-1}$), where $x$ and $y$ are in $M_S$, is 
replaced by a word of length at most 2, the length of the word may decreases or stay the same, this shows 
that the length of the result is not greater than $\ell$. Moreover, since in Lemma 3.9 (ii) of \cite{Deh1}, 
Dehornoy proved that the number of steps needed to $f$-reverse a word of the form $u^{-1}v$ with $u,v \in M_S$ is 
bounded by $\frac{1}{4}N\ell^2$ steps, where $N$ is the upper bound to the number of steps needed to compute 
the $f(x,y)$, we end with a total upper bound for the number of steps which is $2N\ell^2/4=N\ell^2/2$.
\end{proof}

\bt
Activating Algorithm \ref{ALG} takes at most $O(N\ell^2)$.
\et

\begin{proof}
This is an immediate consequence of Proposition \ref{prop4} when noticing that the length of the word before doing 
step $(2)$ of the Algorithm \ref{ALG} may not increase.
\end{proof}

%%%%%%%%%%%%%%%%%%%%%%%%%%%%%%%%%%%%%%%%%%%%%%%%%%%%%%%%%%%%%%%%%%%%%%%%%%%%%%%%%%%%%%%%%%%%%%%%%%%
\section{Application}
We finish the paper with an application for the above algorithm. Recall that braid monodromy is a homomorphism 
between the fundamental group of a punctured disk into the braid group, which is determined by a branch curve 
of a generic projection to $\C ^2$ of a hypersurface in $\C \PP ^n$. It was suggested in \cite{BGTI} and in 
\cite{BGTII} to use the braid monodromy as a tool to compute the fundamental group of the complement of this 
branch curve, and thus to induce an invariant for the surfaces up to diffeomorphism. 

The \emph{braid monodromy factorization} was introduced in \cite{KUL}. It is a product form of the image of the braid monodromy homomorphism on an ordered generating set for fundamental group of the punctured disk.
In \cite{KUL} it is shown that the braid monodromy factorization of the branch curve can serve under 
certain conditions and a specific equivalence relation, as a tool to classify surfaces up to 
diffeomorphism. 
Moreover, in \cite{KUL} [Proposition 3.3] it was proved that if the branch curve is cuspidal, then all 
factors of the braid monodromy factorization must be half-twists to some power. Since the set of all 
half-twists in the braid group is the conjugacy class of its generators one can apply the algorithm presented 
here to each of the factors in the factorization, and test its results. If the result of Algorithm \ref{ALG} 
for at least one of the factors is \emph{False}, then the factorization cannot be a braid monodromy 
factorization of a cuspidal branch curve.

This can be done in $O(N\ell^2 \cdot r)$ where $r$ is the number of elements of the factorization, $\ell$ is 
the length of the longest of all factors and $N$ is as above.

Since in the braid group $M_S$ is actually the set of all positive braids in which no two strings cross more 
than once (\emph{permutation braids}), the number of steps needed to compute $f(x,y)$ or $g(x,y)$ is bounded 
by $O(n\cdot log \ n)$ where $n$ is the braid's number of strings \cite{Word process}. Therefore, in this case we 
have a total of $O(n\cdot log \ n \cdot \ell^2\cdot r)$ operations needed to test the validity of a 
factorization by this method. Finally, using the standard embeddings of the Coxeter groups of type $B(n)$ and 
$D(n)$ in the braid group on $2n+1$ strings (see \cite{Eriksson} for instance), one can replace in Theorem 
\ref{THcomplex} the constant $N$ by $(2n+1) log(2n+1)$ in these cases.

%%%%%%%%%%%%%%%%%%%%%%%%%%%%%%%%%%%%%%%%%%%%%%%%%%%%%%%%%%%%%%%%%%%%%%%%%%%%%%%%%%%%%%%%%%%%%%%%%%%

Godelle Eddy\\
Laboratoire de Math\'ematique Nicolas Oresme,\\
Universit\'e de Caen,\\
14032 Caen cedex, France\\
email: Eddy.Godelle@math.unicaen.fr

\noindent Kaplan Shmuel\\
Department of Mathematics and Statistics,\\
Bar-Ilan University,\\
Ramat-Gan 52900, Israel\\
email: kaplansh@macs.biu.ac.il

\noindent Teicher Mina\\
Department of Mathematics and Statistics,\\
Bar-Ilan University,\\
Ramat-Gan 52900, Israel\\
email: teicher@macs.biu.ac.il\\
\end{document}